\documentclass[12pt,a4paper]{article}
\usepackage{graphicx}
\usepackage{subcaption}
\usepackage{fullpage}

\usepackage[cmex10]{amsmath}
\usepackage{mathtools}
\usepackage{amsthm,amssymb}
\usepackage{amsfonts}
\usepackage{float}
\usepackage[utf8]{inputenc}
\usepackage{makeidx}
\usepackage[english]{babel}
\usepackage{verbatim}
\usepackage[autostyle]{csquotes}
\usepackage[numbers]{natbib}
\usepackage{amsthm}
\usepackage[bottom]{footmisc}
\usepackage{enumitem}
\theoremstyle{plain}

\usepackage{verbatim}
\newtheorem{theorem}{Theorem}%
\newtheorem{definition}{Definition}
\newtheorem{lem}[subsubsection]{Lemma}

\graphicspath{{./figures/}}
\title{Rank-sparsity decomposition for planted quasi clique recovery }
\author{Sakirudeen A. Abdulsalaam and Montaz Ali}
\date{}
\begin{document}
	\maketitle
	\section*{Abstract}
	In this paper, we apply the Rank-Sparsity Matrix Decomposition to the planted Maximum Quasi-Clique Problem (MQCP). This problem has the planted Maximum Clique Problem (MCP) as a special case. The maximum clique problem is NP-hard. A Quasi-clique or $\gamma$-clique is a dense graph with the edge density of at least $\gamma$, $\gamma \in (0, 1]$. The maximum quasi-clique  problem seeks to find such a subgraph with the largest cardinality in a given graph. Our method of choice is the low-rank plus sparse matrix splitting technique.
We present a theoretical basis for when our convex relaxation problem recovers the planted maximum quasi-clique. We have derived a new bound on the norm of the dual matrix that certifies the recovery using $l_{\infty, 2}$ norm. We have showed that when certain conditions are met, our convex formulation recovers the planted quasi-clique exactly. The numerical experiments we have performed corroborate our theoretical findings.
	\section{Introduction}\label{rmpintro}
	Many real-life optimization problems involve constraints on rank of matrices. Low order controller design, minimal realization theory, and model reduction are some examples from control theory \cite{rechtguaranteed}. Furthermore, in some areas of artificial intelligence, the problems of manifold learning, multi-task learning and inference with partial information have been cast as rank minimization problems (RMP) \cite{recht2008necessary}. \citeauthor{amesnuclear} \cite{amesnuclear} recently applied rank minimization technique to solve the maximum clique problem. The maximum clique problem requires finding the largest complete subgraph in a given graph. The maximum clique problem is applicable in identifying connected components in a computer or telecommunications network \cite{abellomaximum}, coding theory, fault diagnosis in multiprocessor systems, and pattern recognition \cite{ames}. However, the maximum clique problem is too confining for some applications like data mining/clustering, community detection, criminal network analysis, epidemic control \cite{balasthesis}. Therefore, various clique relaxation models have emerged to bridge this gap. Prominent among them is the maximum quasi-clique  also known as maximum $\gamma$-clique. Finding maximum $\gamma$-clique in a graph means finding the largest subgraph with required edge density, $\gamma \in (0, 1]$. When $\gamma = 1$, the maximum quasi-clique problem is equivalent to the maximum clique problem. Various heuristic methods are available for maximum quasi-clique recovery but the theoretical basis is lacking. This paper presents the first theoretical framework for maximum quasi-clique recovery.  We formulate a new convex problem for $\gamma$-clique recovery by generalising the idea in \cite{amesnuclear}. We formulate a convex model for recovery of planted maximum quasi-clique using the nuclear norm and matrix $l_1$ norm (defined in Section \ref{notations}). More importantly, we establish a guarantee of convergence for the formulated convex program through the application of a novel norm, the \(l_{\infty,2}\) norm. Finally, we have applied the formulated model in recovery of planted maximum quasi-clique using a sythetic data.
	
The rest of this paper is organised as follows. We define symbols and notations in Section \ref{notations} while a brief survey of previous related work is contained in Section \ref{lr}. We present our proposed convex program and the proof of our main result in Section \ref{mainresults} and \ref{mainproof} respectively. Section \ref{experiments} contains the numerical experiment while we give concluding remarks in Section \ref{conclusion}.

\section{Notations and Definitions}\label{notations}
	Before we proceed, we summarize the notations used in this article. We use capital letter $X$ to denote a matrix $X \in \mathbb{R}^{n_1 \times n_2}$. $X^T$ is the transpose of $X$. The singular values of $X$ is the square root of the eigenvalues of the matrix $XX^T$. If $X$ is positive semidefinite, then the singular values are equivalent to the eigenvalues of $X$. $\sigma _i(X)$ denotes the $i$-th largest singular value of $X$.  For a $d$-dimensional vector $v \in \mathbb{R}^d$, its Euclidean norm is defined as $||v|| = \left(\sum_{i = 1}^{d} {v_i}^2 \right)^{1/2}$. Let $X$ and $Y$ be matrices with the same dimension, the Euclidean inner product between $X$ and $Y$ is defined as $\langle X,Y\rangle := trace(X^TY) = \sum_{i = 1}^{n_1} \sum_{j = 1}^{n_2} X_{ij}Y_{ij}.$ This inner product is associated with a norm known as the Frobenius norm or Hilbert-Schmidt norm and is defined as $||X||_F := \langle X , X \rangle ^{1/2} = \left(\sum_{i = 1}^{n_1} \sum_{j = 1}^{n_2} X_{ij}^2\right)^{1/2} = \left( \sum _{i = 1} ^ r \sigma_i ^2 \right) ^{1/2},$ where $r$ is the rank of $X$ (the number of linearly independent rows or columns of $X$ or the number of its non-zero singular values). The sum of the singular values of a matrix is known as the nuclear norm of the matrix and is defined as 
	$||X||_* := \sum_{i = 1}^{r} \sigma_i(X).$ The norm is also known by some other names which include the schattern  $1$-norm, the Ky Fan $r$-norm and the trace class norm \cite{rechtguaranteed}. The spectral norm, operator norm or induced $2$-norm of a matrix is its largest singular value. It is defined as $||X|| := \sigma_1(X).$ The three norms defined above are related by the following inequalities \cite{rechtguaranteed}:
	\begin{equation} 
	||X|| \leq ||X||_F \leq ||X||_* \leq \sqrt{r}||X||_F \leq r ||X||. \label{normsinequality}
	\end{equation}

	From equation \eqref{normsinequality}, we have that
	$$rank(X) \geq \frac{||X||_*}{||X||}, \quad \forall \text{ } X.$$
	The last two matrix norms that we will define here are the $l_1$ and the $l_\infty$ norm. For a matrix $X$, $||X||_1 := \sum_{i j} |X_{ij}|$ while $||X||_\infty = \max_{ij} |X_{ij}|$.
	
	\section{Related works}\label{lr}
	We give summary of the existing related works in the section. 
The heuristic based on nuclear norm minimization has been introduced in \cite{rechtguaranteed}. It is an extension of the well known trace heuristic (see for example \cite{beckcomputational,mesbahirank}), used in control theory, for rank minimization. 
The affine rank minimimization problem can be expressed as:
	\begin{equation}\label{rankmodel}
	\begin{aligned} 
	&\text{minimize } rank(X),\\ 
	&\text{subject to } \mathcal{A}(X) = b, 
	\end{aligned}
	\end{equation}
	where $X \in \mathbb{R}^{n_1 \times n_2}, b \in \mathbb{R}^d$ and the linear map $\mathcal{A} : \mathbb{R}^{n_1 \times n_2} \mapsto \mathbb{R}^d$.
	The rank minimization problem is NP-Hard and non-convex in general. All the available finite time algorithms have exponential running time. Nevertheless, in some cases with special structure, RMP can be reduced to solution of a linear system or solved by the singular value decomposition approach \cite{rechtguaranteed}. 
	When the singular values are all equal to one, then $rank(X)$ is equal to the nuclear norm, i.e, $rank(X) = ||X||_*$. When the singular values are less than or equal to one, then $rank(X) \geq ||X||_*$. Hence, the nuclear norm is the tightest convex function that underestimates the rank function on the unit ball in the spectral norm. 
	Therefore, the nuclear norm heuristic provides bounds for solutions of affine rank minimization problem \cite{rechtguaranteed}. 
	
	
	Hence, the objective function \eqref{rankmodel} can be replaced with nuclear norm and then we have the following formulation:
	\begin{equation}\label{nnmodel}
	\begin{aligned} 
	&\text{minimize } ||X||_*, \\
	&\text{subject to } \mathcal{A}(X) = b. 
	\end{aligned}
	\end{equation}
Indeed, when $X$ is positive semidefinite, \eqref{nnmodel} is equivalent to trace heuristic \cite{beckcomputational,mesbahirank}. The convex problem \eqref{nnmodel} can be solved using a wide range of existing algorithms. Over time, nuclear norm has been observed to produce very low rank solutions in practice but the theoretical basis for when it does produce the minimum rank solution recently emerged \cite{rechtsimpler}. 

An instance of constrained matrix rank minimization problem is the matrix completion problem: recovering $M \in \mathbb{R}^{n_1 \times n_2}$ of rank $r$ from its $m$ observed entries. Recovery is possible for $m > r(n_1 + n_2 - r)$ \cite{candesexactmatrix,candespower}. In addition, when $M$ is of low-rank, the following incoherence conditions \cite{candesrobust,chandrasekaranrank,grossrecovering} (which ensures that the column/row spaces of $M$ are not closely aligned with the cannonical basis vectors ${e_i}$), with parameter $\mu \in [1, \frac{\min(n_1, n_2)}{r}]$, hold: 

\begin{equation}\label{standardincoherence}
	\max_i ||U^T e_i||^2 \leq \frac{\mu r}{n_1}, \quad \max_i ||V^T e_i||^2 \leq \frac{\mu r}{n_2}
	\end{equation}
	and
	\begin{equation}\label{jointincoherence}
	||UV^T||_\infty \leq \sqrt{\frac{\mu r}{n_1 n_2}},
	\end{equation} 
where $M = \sum_{i = 1}^{r} \sigma_i u_i v_i^T = U\Sigma V^T$ and $\sigma_1, \ldots, \sigma_r$ are singular values of $M$.	
	
	An extension of the matrix completion problem  known as matrix decomposition problem was studied in \cite{candesrobust,chandrasekaranrank,chenlow}.
In this setting, the low rank matrix is a submatrix of a larger matrix with other additional unwanted entries (sometimes called error or corruptions). Such a matrix is formed by adding a low-rank matrix to a sparse matrix. The goal is to recover the low-rank and sparse components of the given matrix. A more difficult version of this is the case when the low-rank matrix has some missing entries. \citeauthor{chenlow} \cite{chenlow}  calls this case recovery in presence of errors and erasures. Generally, the matrix decomposition problem is stated as follows. We are given a matrix $M \in \mathbb{R}^{n_1 \times n_2}$, which is a sum of a low rank matrix $B \in \mathbb{R}^{n_1 \times n_2}$ and a sparse (errors) matrix $C \in \mathbb{R}^{n_1 \times n_2}$. The cardinality and the values of the non-zero entries of $C$ are not known; likewise the locations of the non-zero entries are also not known. The goal is to recover $B$ and $C$ from $M$. The convex formulation for this problem is:
\begin{equation}\label{mdp}
\begin{aligned}
\min \text{ } & ||B||_* + \lambda ||C||_1, \\
\text{subject to } & B + C = M,
\end{aligned}
\end{equation}
where $\lambda$ is a regularisation parameter. We will discuss the choice of its value in the next section. 
It turns out that in addition to the incoherence conditions \eqref{standardincoherence} -- \eqref{jointincoherence}, it is also required that  sparsity pattern does not exist in the entries of $C$. Moreover, exact recovery of $(B, C)$ via \eqref{mdp} is theoretically linked with finding an approximate dual matrix, $Q$, for which the Golfing Scheme is used \cite{candesrobust}.

 In this paper, we adapt matrix decomposition technique to recovery of planted quasi-clique. We are inspired by the work of \citeauthor{amesnuclear} \cite{amesnuclear} but our problem formulation and proof technique differ greatly. We provide a novel theoretical guarantee for maximum quasi-clique recovery using convex relaxation. The major tool used to achieve the new results was the $l_{\infty, 2}$ norm. 
Our work can be seen as an extension of \cite{chenincoherence} to matrix decomposition problem.
	
\section{Mathematical model and recovery of planted quasi-clique}\label{mainresults}
The maximum clique problem (MCP) is a combinatorial problem with applications in various fields \cite{abellomaximum, amesconvex, ostergard}. An instance of MCP is the planted (hidden) clique. One of the ways to plant a clique is the randomized approach where a graph of size $n$ is constructed and a clique of size $n_c$ is inserted. The remaining edges are then constructed with probability $\rho$.

Inadequacies of the clique model for some applications led to the emergence of version of clique relaxations. These are $k$-clique \cite{luceconnectivity}, $k$-club \cite{albagraph,mokkencliques} and $k$-clan \cite{mokkencliques}; all are reachability (diameter) based relaxations. On the other hand, $k$-plex \cite{seidmangraph,balascliquerelaxations,balasthesis} is a degree based relaxation. The case $k = 1$ in this model is equivalent to the classical clique model.
 In this paper, we focus on a relatively new clique relaxation model known as \emph{quasi-clique} or $\gamma$-clique; $\gamma \in (0, 1]$ is the relaxation parameter. When $\gamma = 1$, the quasi-clique model is equivalent to clique model. Therefore, the quasi-clique model generalizes the clique model. Quasi-clique is an edge density based relaxation. It was proposed by \citeauthor{abellomaximum} \citep{abellomaximum, abellomassive}.  For a subgraph $C_\gamma = (V', E')$ of $G = (V, E)$ to be a $\gamma$-clique, the inequality $|E'| \geq \gamma \binom{n_c}{2}$, where $n_c = |V'|$, must be satisfied. Furthermore, for the inequality to hold, the degree of every $v \in V'$ must be at least $\lceil \frac{\gamma(|V| - 1)}{2} \rceil$. Consequently, $\gamma$-clique provides a mix of relaxation in reachability, degree and robustness. As a result, the model fits more real life situations while preserving the ideal properties of a clique as much as possible.  

Majority of the previous studies on $\gamma$-clique focused on developing heuristic based algorithm for detection of large quasi-clique \cite{pattillomax}. \citeauthor{abellomaximum} \citep{abellomaximum} was the first to publish on maximum quasi-clique problem (MQCP). They proposed a greedy randomized adaptive search procedure (GRASP) for finding large quasi-clique in graphs generated from communication data. \cite{abellomassive} presents similar approach using semi-external memory algorithms that handles massive graph together with GRASP so as to be able to deal with graphs with millions of nodes. \cite{brunatoeffective} presents a heuristic solution approach to quasi-clique problem by extending two stochastic local search algorithms for the classical MCP to $\gamma$-clique problem. 
The first study on quasi-clique from mathematical perspective was carried out in \cite{pattillomax} where $\gamma$-clique problem was proved to be NP-complete. Also, an upper bound was derived for maximum $\gamma$-clique and mixed integer programming models were proposed for finding maximum $\gamma$-clique in a random graph. Recently, the model in \cite{pattillomax} has been modified and improved in \cite{veremyevexact}.
In both cases, only numerical experiments were presented with no theoretical guarantee for perfect quasi-clique recovery. 

\citeauthor{amesnuclear} \cite{amesnuclear} proposed the nuclear norm relaxation of the planted clique problem:

\begin{subequations}\label{mcpnn}
	\begin{align}\label{mcpnn:a}
	&\min ||X||_*, \\
	\label{mcpnn:b}\text{subject to } & \sum_{i \in V} \sum_{j \in V} X_{ij} \geq n_c^2, \\
	\label{mcpnn:c} & X_{ij} = 0 \text{ } \forall \text{ } ij \notin E \text{ and } i \neq j,
	\end{align}
\end{subequations}
where $G = (V,E)$. The sufficient conditions of recovery provided for this model are based on the size of the graph, the size of the planted clique and the probability of the diversionary edges. Under these conditions, recovery is guaranteed provided the hidden clique is large enough.

We propose the following relaxed problem for the planted $\gamma$-clique problem:
\begin{subequations}\label{qcp}
\begin{align}\label{qcp:a}
&\min ||B||_* + \lambda ||C||_1 \\
\label{qcp:b}\text{subject to } & \sum_i \sum_j B_{ij} \geq \gamma n_c ^2   \\
\label{qcp:d}& B + C = M \\ 
\label{qcp:f}& B_{ij}, C_{ij} \in [0, 1], \quad n_c \in \mathbb{N}, 	
\end{align}
\end{subequations}
where $M$ is the adjacency matrix of the given graph $G$; $B$, $C$ and $n_c$ are optimization variables. Our formulation is related to matrix decomposition. The idea is to split $M$ into $B$ (the adjacency matrix corresponding to the planted quasi-clique) and a sparse matrix $C$ (those non-zero entries corresponding to the edges not in the quasi-clique). The constraint \eqref{qcp:b} ensures that the recovered quasi-clique satisfies the edge density requirement while 
\eqref{qcp:d} makes sure that the decomposition agrees with the input matrix. $n_c$ is an integer value variable that determines the size of the recovered quasi-clique.

Since we are only interested in $B$, we can eliminate constraint \eqref{qcp:d} and write $C = M - B$. Therefore, \eqref{qcp} can be reformulated as

\begin{subequations}\label{qcp1}
	\begin{align}\label{qcp1:a}
	&\min ||B||_* + \lambda ||M - B||_1 \\
	\label{qcp1:b}\text{subject to } & \sum_i \sum_j B_{ij} \geq \gamma n_c^2   \\
	\label{qcp1:e}&B_{ij} \in [0, 1], \quad n_c \in \mathbb{N}. 	
	\end{align}
\end{subequations}

Clearly, our model differs from that of \citeauthor{amesnuclear} \cite{amesnuclear}. A striking difference between problem \eqref{mcpnn} and \eqref{qcp1} is that \eqref{mcpnn} requires $n_c$ as an input because the size of the clique is required while we treat $n_c$ as an optimization variable. Hence, the proof technique for our recovery guarantee is different. Our approach is based on matrix decomposition problem. The only point of convergence between the two methods is the prerequisite that recovery is feasible only if $B$ is of low rank. We propose the following theorem.
\begin{theorem}\label{main}
Let $Q_\gamma = (V', E')$ be an $n_c$-vertex $\gamma$-clique planted in a graph $G = (V, E)$ of $n$-vertices ($n_c < n$) with the adjacency matrix $M = B + C$. Suppose $B \in \mathbb{R}^{n \times n}$, the adjacency matrix of $Q_\gamma$, obeys the incoherence condition. Let $\lambda = \frac{1}{\sqrt{n}}$ and denote the optimal pair of problem \eqref{qcp} as $(B^*, C^*)$. Then there exists universal positive constants $c$ and $c_0$; independent of $n$, such that if
$$c_0 \frac{\mu r \log n}{n} \leq p \leq 1,$$
where $p$ denotes the sampling probability of entries of  $B$, then $(B^*,C^*)$ is the
 unique optimal solution of \eqref{qcp} with probability at least $1 - c_0n^{-10}$.
\end{theorem}

\section{Proof of Theorem  \ref{main}}\label{mainproof}
	In this section, we present the proof of our main result, Theorem \ref{main}. The main steps follow the general idea in the low-rank matrix recovery literature \cite{candescomletion_with_noise, candesexactmatrix, candesrobust, chenlow, licompressed}. The steps involved include constructing a dual matrix $Q$, which certifies the optimality of $(B^*, C^*)$ for the convex problem (\ref{mdp}). This dual certificate must obey some subgradient-type conditions. One of these conditions is that the spectral norm of $Q$, i.e $||Q||$, must be small. In many of the previous approaches, $||Q||$ is bounded by the matrix $l_\infty$ norm. Since $Q$ has singular value decomposition, $Q = U \Sigma V$, by virtue of equation \eqref{jointincoherence}, there is a relationship between the bound on the dual certificate and the incoherence condition. We bound the norm of the dual certificate by the $l_\infty$ and $l_{\infty, 2}$ norm. The $l_{\infty, 2}$ was first used in \cite{chenincoherence} to derive a tighter bound for matrix completion problem. Ours is an extension of this concept to the matrix decomposition setting. The $l_{\infty, 2}$ norm is defined on a matrix $M$ as:
	\begin{equation}
	||M||_{\infty, 2} := \max \left\{\max_i \sqrt{\sum_{b} M_{ib}^2}, \max_j \sqrt{\sum_{a} M_{aj}^2} \right\}. \label{elinf}
	\end{equation}
	It is noteworthy that for any matrix $M \in \mathbb{R}^{n_1 \times n_2}, ||M||_{\infty, 2} \leq \sqrt{\max \{n_1, n_2\}} ||M||_\infty$ \cite{chenincoherence}. Therefore, from the incoherence property, we have
	\begin{equation*}
	||UV^T||_{\infty, 2} \leq \sqrt{\max \{n_1, n_2\}} ||UV^T||_\infty \leq \sqrt{\frac{\mu r}{\min\{n_1, n_2\}}},
	\end{equation*}
	or
	\begin{equation}\label{incoherenceforlinfinity2norm}
	||UU^T||_{\infty, 2} \leq \sqrt{n}||UU^T||_\infty \leq \sqrt{\frac{\mu r}{n}},
	\end{equation}
	for a square matrix.
	
	There is a norm similar to the $l_{\infty, 2}$ norm which is the $l_{2, \infty}$ norm defined as $||M||_{2, \infty} := \max_{x \neq 0}\frac{||M x||_\infty}{||x||_2}.$ 
	The $l_{2, \infty}$ is the maximum Euclidean norm of the rows of $M$ while the $l_{\infty, 2}$ is the maximum of both the rows and column norm of $M$ \cite{chenincoherence, rebrovanorms}. These norms yield a tighter bound on the entries of a matrix than the commonly used ones \cite{capetwo, chenincoherence}. We show that the $l_{\infty, 2}$ norm yields a tighter bound on the norm of the dual matrix. This is achieved by expressing the bound as a sum of the $l_\infty$ and $l_{\infty, 2}$ norm instead of the previously derived bounds in matrix decomposition literature which use the $l_\infty$ norm only. This is one of the main contributions of this paper. In addition, we adopt some novel simplifying ideas different from the existing works.

	Since every graph has a square adjacency matrix, we prove only for square matrices. The arguments follow easily for rectangular case by taking $n = \min{\{n_1, n_2\}}$. We denote universal constants that do not depend on the parameters of the problem (i.e $n, r, \mu$, etc) with $c_0, c_1 \text{ and } c_2$. By with high probability, we mean with probability at least $1 - c_1n^{-c_2}, c_1, c_2 > 0$. 
Suppose $B \in \mathbb{R}^{n \times n}$ is a symmetric matrix of rank $r$. Then, $B$ is orthogonally diagonalizable. Hence $B = U\Sigma U^T; U = (u_1,u_2, \ldots, u_r)$ and $\Sigma = diag(\sigma_1, \ldots, \sigma_r)$, where $u_i$ and $\sigma_i$ are the $i$-th singular vector and singular value, respectifully. The subgradient of the nuclear norm at $B$ is of the form:
$$UU^T +Q,$$
where $U^TQ = 0, QU = 0 \text{ and } ||Q|| \leq 1.$

In addition, we define $T$, the linear space of matrices that share the same row space or column space as $B$, as 
\begin{equation}
T := \{ UX^T + YU^T : X,Y \in \mathbb{R}^{n \times r} \},
\end{equation}
It can easily be shown that $\mathcal{P}_{T^\perp}Q = Q$. It follows that, for any matrix $Z \in \mathbb{R}^{n \times n}$
\begin{equation}\label{equ:mat_ort_pro}
\mathcal{P}_{T^\bot}Z = (I - UU^T)Z(I - UU^T),
\end{equation}
 where $I -UU^T$ is the orthogonal projections onto the orthogonal complement of the linear space spanned by the columns of $U$. As a result of this, for any matrix $Z$, $||\mathcal{P}_{T^\bot}Z|| \leq ||Z||.$ We will make use of this fact later.
 
The subgradient of $l_1$-norm at $C$, where $\Gamma = supp(C)$, is of the form:
$$Sgn(C) + D$$
where $\mathcal{P}_\Gamma D = 0$ and $||D||_\infty \leq 1$. Here, support of $C$ is defined by $supp(C) = \{ (i,j) | C_{i,j} \neq 0\}$ and $Sgn(.)$ is the sign function.

\subsection{The Bernoulli Model and Derandomization}
\subsubsection{Bernoulli Model}
Rather than showing that Theorem \ref{main} holds with $\Gamma$ sampled uniformly, where $\Gamma$ is a random subset 
of cardinality $k$, it is easier to prove the theorem for $\Gamma$ sampled according to the Bernoulli model, with $\Gamma = \{(i, j): \Delta_{ij} = 1 \}$, where $\Delta_{ij}$'s are independent and identically distributed Bernoulli random variables. $\Delta_{ij} = 1$ with probability $\varrho$ and $0$ with probability $1 - \varrho$. Hence, $\varrho n^2$ is the expected cardinality  of $\Gamma$. By doing this, we exploit the statistical independence of measurements. Based on the arguments presented in \cite{candescomletion_with_noise,candesexactmatrix, candesrobust}, any theoretical guarantee proved for Bernoulli model is also valid for the uniform model and the converse holds as well. Henceforth, we will write $\Gamma \sim Bern(\varrho)$, i.e, $\Gamma$ follows Bernoulli distribution with parameter $\varrho$.

\subsubsection{Derandomization}
Note that the sign of the entries of $C$ in Theorem \ref{main}  is fixed. Surprisingly, it is more convenient to prove the theorem under a tougher condition. We assume that the sign of ${C}_{ij}$, for ${C}_{ij} \neq 0$, are independent symmetric Bernoulli random variables which assume the value $1$ or $-1$ with probability $1/2$. The probability of recovering $C$ with random sign on the support set of $\Gamma$ is at least the same with $C$ with fixed sign. This randomization technique was invented in \cite{candesrobust} and has been previously used in \cite{nguyenexact} and \cite{licompressed}.  The supporting theorem is stated formally here.

\begin{lem}[Theorem 2.3 of \cite{candesrobust}]
Suppose $B$ obeys the conditions of Theorem \ref{main} and that the positions of the non-zero entries of $C$ follow the Bernoulli model with parameter $2\varrho$, with the signs of $C$ independent and identically distributed $\pm 1$ with probability $1/2$. Then if the solution to problem (\ref{qcp}) is exact with high probability, it is also exact for the model with fixed signs and location sampled from the Bernoulli model with parameter $\varrho$ with at least the same probability.
\end{lem}

\subsection{Subgradient Condition for Optimality}
Now, we state sufficient conditions for the pair $(B^*, C^*)$ to be the unique optimal solution to problem (\ref{mdp}). The conditions are stated in terms of the dual matrix, whose existence certifies optimality. These conditions are given in the following lemma which is similar to Proposition $2$ of \cite{chandrasekaranrank} and Lemma $2.4$ of \cite{candesrobust}.

\begin{lem}[Proposition 2 of \cite{chandrasekaranrank}, Lemma 2.4 of \cite{candesrobust}]\label{lemma2.4candes}
	Let $M = B + C$. Suppose $\Gamma \cap T = \{ 0 \}$ (i.e, $||\mathcal{P}_\Gamma \mathcal{P}_T|| < 1)$. Then $(B^*, C^*)$ is the unique optimal solution to (\ref{qcp}) if there exists a pair of dual matrix $(D, Q)$ such that
	$$UU^T + Q = \lambda (Sgn(C) + D),$$
	with $\mathcal{P}_T Q = 0, ||Q|| < 1, \mathcal{P}_\Gamma D = 0 \text{ and } ||D||_\infty < 1.$
\end{lem}

The proof of Lemma \ref{lemma2.4candes} can be found in Appendix \ref{app1}. Consequently, to prove exact recovery, it suffices to derive a dual certificate $Q$ which satisfies:
\begin{enumerate}
	\item[(a)] $Q \in T^\perp$,
	\item[(b)] $||Q|| < 1$,
	\item[(c)] $\mathcal{P}_\Gamma (UU^T + Q) = \lambda Sgn(C)$,
	\item[(d)] $||\mathcal{P}_{\Gamma^\perp}(UU^T + Q)||_\infty < \lambda$.
\end{enumerate}
Our approach, however, yields a somewhat different certificate with high probability. This is due to a relaxation on the constraint $\mathcal{P}_{\Gamma}(UU^T + Q) = \lambda Sgn(C)$. This relaxation was introduced in \cite{grossrecovering} and was previously used in \cite{candesrobust}.

\begin{lem}[Lemma 2.5 of \cite{candesrobust}]\label{candeslemma25}
Suppose $||\mathcal{P}_\Gamma \mathcal{P}_T|| \leq 1/2$ and $\lambda < 1$. Then $(B^*, C^*)$ is the unique optimal solution of (\ref{qcp}) if there exists a pair of dual matrix $(Q, D)$, such that
$$UU^T + Q = \lambda (Sgn(C) + F + \mathcal{P}_\Gamma D)$$
with $\mathcal{P}_T Q = 0, \text{ } ||Q|| \leq 1/2, \text{ } \mathcal{P}_\Gamma D = 0, \text{ } ||F||_\infty \leq 1/2, \text{ and } ||\mathcal{P}_\Gamma F|| \leq 1/4.$
\end{lem}

Interested reader can check the proof of this Lemma in \cite{candesrobust}. As a result of Lemma \ref{candeslemma25}, it is sufficient to derive a dual matrix $Q$, which obeys
\begin{equation}\label{stasta}
\begin{aligned}
&Q \in T^\perp,\\
& ||Q|| < 1/2,\\
& ||\mathcal{P}_\Gamma(UU^T - \lambda Sgn(C) + Q)||_F \leq \lambda / 4, \\
&||\mathcal{P}_{\Gamma^\perp}(UU^T + Q)||_\infty < \lambda/2. 
\end{aligned} 	
\end{equation}
\subsection{Construction of Dual Certificate}\label{sec:construction}
We need to construct a dual matrix $Q$ satisfying the conditions in (\ref{stasta}). $Q$ will be constructed using a modified version of the \emph{Golfing Scheme} \cite{candesrobust, chenincoherence, chenlow, grossrecovering}. Fix the value of $k_0$, e.g set $k_0 = 20 \lceil \log n \rceil$ and decompose the observed entries of $\Gamma$ into independent entries. Recall that $\Gamma \sim Bern(p)$.  Therefore, $\Gamma ^ C \sim Bern(1 - p)$. $\Gamma ^ C$ is a union of $k_0$ independent samples, i.e, $\Gamma^C = \bigcup _{1 \leq k \leq k_0} \Gamma_k$, with each $\Gamma_k$ following the Bernoulli model with parameter $q$ (probability of sampling for each batch). This gives the following Binomial model:
\begin{equation}\label{eq:binombatch}
\mathbb{P} \left( (i, j) \in \Gamma \right)  = \mathbb{P} ( Bin (k_0, q) = 0) = (1 - q) ^ {k_0},
\end{equation}
where $Bin$ in equation \eqref{eq:binombatch} means binomial. The two models are equivalent if $p = (1 - q)^{k_0}$. 

We now proceed to construct a dual certificate $Q = Q_B + Q_C$. Details about the two components of $Q$ are as follow.

\subsubsection{Constructing $Q_B$ using Golfing Scheme}
Let $k_0$ and $\Gamma _k, k \in [1, k_0]$, be as defined above. Also, recall that $\Gamma ^ C = \cup _{1 \leq k \leq k_0}\Gamma_k$. Starting with $Y_0 = 0$ and proceeding inductively, we define
\begin{equation}\label{golfingscheme}
Y_k = Y_{k - 1} + p^{-1}\mathcal{P}_{\Gamma_k} \mathcal{P}_T(UU^T - Y _{k - 1}),
\end{equation}

and set
\begin{equation}\label{eq:qbequalyk0}
Q_B = \mathcal{P}_{T^\perp}Y_{k_0}.
\end{equation}

\subsubsection{Constructing $Q_C$ using Least Squares Method}
Suppose $||\mathcal{P}_\Gamma \mathcal{P}_T|| < 1/2$, then $||\mathcal{P}_\Gamma \mathcal{P}_T \mathcal{P}_\Gamma|| < 1/4$. Thus, the operator $\mathcal{P}_\Gamma - \mathcal{P}_\Gamma \mathcal{P}_T \mathcal{P}_\Gamma$ is invertible.  We write $(\mathcal{P}_\Gamma - \mathcal{P}_\Gamma \mathcal{P}_T \mathcal{P}_\Gamma)^{-1}$ for the inverse and set
\begin{equation}
Q_C = \lambda \mathcal{P}_{T^\perp}(\mathcal{P}_\Gamma - \mathcal{P}_\Gamma \mathcal{P}_T \mathcal{P}_\Gamma)^{-1} Sgn(C), \label{defq_c}
\end{equation}
according to the least squares method  \cite{candesrobust}.
There is an alternative definition of (\ref{defq_c}) using the convergent Neumann series \cite{candespower, candesrobust}. The definition is as follows:
\begin{equation}
Q_C = \lambda \mathcal{P}_{T^\perp} \sum _{k \geq 0}(\mathcal{P}_\Gamma \mathcal{P}_T \mathcal{P}_\Gamma)^k Sgn (C). \label{defq_c_alternative}
\end{equation}
Observe that
\begin{align}
\mathcal{P}_\Gamma Q _ C &= \lambda \mathcal{P}_\Gamma \mathcal{P}_{T^\perp}(\mathcal{P}_\Gamma - \mathcal{P}_\Gamma \mathcal{P}_T \mathcal{P}_\Gamma)^{-1} Sgn(C), \nonumber \\
&= \lambda \mathcal{P}_\Gamma(\mathcal{I} - \mathcal{P}_T)(\mathcal{P}_\Gamma - \mathcal{P}_\Gamma \mathcal{P}_T \mathcal{P}_\Gamma)^{-1} Sgn(C), \nonumber \\
&= \lambda (\mathcal{P}_\Gamma - \mathcal{P}_\Gamma \mathcal{P}_T \mathcal{P}_\Gamma)(\mathcal{P}_\Gamma - \mathcal{P}_\Gamma \mathcal{P}_T \mathcal{P}_\Gamma)^{-1} Sgn(C),\label{heart}\\
&= \lambda Sgn(C). \nonumber
\end{align}
\eqref{heart} follows from properties of the projection operators, i.e for any two orthogonal projection operators $\mathcal{P}_1$ and $\mathcal{P}_2$, $\mathcal{P}_1 = \mathcal{P}_1^2$ while $\mathcal{P}_1\mathcal{P}_2 = \mathcal{P}_2\mathcal{P}_1$ must hold for the product $\mathcal{P}_1 \mathcal{P}_2$ to be a projector \cite{baksalarycommutativity}.
Consequently, one can verify that among all the matrices $Q \in T^\perp$ which satisfies $\mathcal{P}_\Gamma Q = \lambda Sgn(C)$, $Q_C$ is the one with minimum Frobenius norm \cite{candesrobust}. Since, by construction, $Q_B, Q_C \in T^\perp$ and $\mathcal{P}_\Gamma Q_C = \lambda Sgn(C)$, it remains to show that $Q = Q_B + Q_C$ obeys:
\begin{align*}
&||Q_B + Q_C|| < 1/2, \\
&||\mathcal{P}_{\Gamma}(UU^T + Q_B)||_F \leq \lambda / 4, \\
& ||\mathcal{P}_{\Gamma^\perp}(UU^T + Q_B + Q_C)||_\infty < \lambda / 2,
\end{align*}
to establish \eqref{stasta}. Indeed, our approach yields a better bound on $Q$. Our new results is proposed in the following lemma.
\begin{lem}\label{mymainlemma}
	Suppose $\Gamma \sim Bern(p)$ such that $0 < p \leq 1$. Let $k_0 = 20 \lceil \log n \rceil$, then under the assumption of Theorem \ref{main}, the matrix $Q_B$ satisfies:
	\begin{enumerate}
		\item[i.] $||Q_B|| < 1/8$,
		\item[ii.] $||\mathcal{P}_\Gamma(UU^T + Q_B)||_\infty < \lambda /8$,
		\item[iii.] $||\mathcal{P}_{\Gamma^\perp}(UU^T + Q_B)||_\infty < \lambda / 4.$
	\end{enumerate}
	Furthermore, if $supp(C) = \Gamma$ ($\Gamma$ is as sampled earlier at the beginning of Section \ref{sec:construction}), and the signs of $C$ are symmetric \emph{iid}, then when the assumptions of Theorem \ref{main} holds, the matrix $Q_C$ satisfies:
	\begin{enumerate}[resume]
		\item[iv.] $||Q_C|| < 1/8$,
		\item[v.] $||\mathcal{P}_{\Gamma^\perp}Q_C||_\infty < 1/4.$
	\end{enumerate}
\end{lem}
\subsection{Key Lemmas}
We now state some important Lemmas which are used to establish our results, Theorem \ref{main}. Our proof differs here, substantially, from the existing works. We derive bound on $Q$  in terms of $l_{\infty, 2}$ norm. This is done with the aid of the following two Lemmas. The first lemma is used to bound the operator norm of $(p^{-1} \mathcal{P}_\Gamma - \mathcal{I})Z$, for a matrix $Z$, in terms of the $l_{\infty, 2}$ and $l_\infty$ norms of $Z$. The bound produced this way is tighter than the previous ones \cite{candesrobust,chenlow}.

\begin{lem}[Lemma 2 of \cite{chenincoherence}]\label{lemma2yudongchen}
	Let $Z$ be a square matrix of dimension $n$. For a universal constant $c > 1$, we have
	\begin{equation*}
	||(p^{-1}\mathcal{P}_\Gamma - \mathcal{I})Z|| \leq c \left(\frac{\log n}{p}||Z||_\infty + \sqrt{\frac{\log n}{p}}||Z||_{\infty, 2} \right),
	\end{equation*}
	with high probability.
\end{lem} 
The next Lemma provides for further controls on the $l_{\infty, 2}$ norm.

\begin{lem}[Lemma 3 of \cite{chenincoherence}]\label{lemma3chencoherence}
	Suppose $Z \in \mathbb{R}^{n \times n}$ is a fixed matrix. If $p \geq c_0 \frac{\mu r \log n}{n}$, for some $c_0 > 0$ which is large enough, then
	\begin{equation*}
	||(\mathcal{P}_\Gamma - p^{-1}\mathcal{P}_T \mathcal{P}_\Gamma)Z||_{\infty, 2} \leq 1/2 \sqrt{\frac{n}{\mu r}}||Z||_\infty + 1/2 ||Z||_{\infty, 2},
	\end{equation*}
	with high probability.
\end{lem}
The proof of these two lemmas can be found in \cite{chenincoherence}. We also need the following standard results  which enable  us to manipulate the $l_\infty$ norm.

\begin{lem}[Lemma 4 of \cite{chenincoherence}, Lemma 3.1 of \cite{candesrobust}, Lemma 13 of \cite{chenlow}] \label{lemma4chenincoherence}
	Suppose $Z \in \mathbb{R}^{n \times n}$ is a fixed matrix and that $\Gamma \sim Bern(p)$. If $p \geq c_0 \frac{\mu r \log n}{n}$, for some $c_0 > 0$ large enough, then
	
	\begin{equation*}
	||(\mathcal{P}_T - p ^{-1}\mathcal{P}_T \mathcal{P}_\Gamma \mathcal{P}_T)Z||_\infty \leq 1/2 ||Z||_\infty,
	\end{equation*}
	with high probability.
\end{lem}
\begin{lem}[Theorem 6.3 of \cite{candesexactmatrix}, Lemma 3.2 of \cite{candesrobust}]
	Suppose $Z \in \mathbb{R}^{n \times n}$ is a fixed matrix and that $\Gamma \sim Bern(p)$. If $p \geq c_0 \frac{\mu r \log n}{n}$, for some $c_0 > 0$ sufficiently large, then with high probability,
	\begin{equation*}
	||(p^{-1}\mathcal{P}_{\Gamma } -  \mathcal{I})Z|| \leq c_0' \sqrt{\frac{n \log n}{p}}||Z||_\infty,
	\end{equation*} 
	for some small numerical constant $c_0' > 0$.
\end{lem}

\begin{lem}[Theorem 4.1 of \cite{candesexactmatrix}, Theorem 2.6 of \cite{candesrobust}, Lemma 1 of \cite{chenincoherence}]\label{lemma1chenchoerence}
	Suppose $p \geq c_0 \frac{\mu r \log n}{n}$ for some sufficiently large $c_0$, then (with high probability);
	$$||\mathcal{P}_T - p ^{-1}\mathcal{P}_T \mathcal{P}_\Gamma \mathcal{P}_T|| \leq 1/2.$$
\end{lem}
We will also need the following results on bounds of the operator norm of random matrices.

\begin{lem}[Corollary 2.3.5 of \cite{taotopics}]\label{terry235}
	Let $Z \in \mathbb{R}^{n \times n}$ be a random matrix with $Z_{ij}$ being independent and identically distributed random variables. Suppose $Z_{ij}$ is uniformly bounded in magnitude by one and has mean zero. Then there exist constants $c, c^* > 0$ such that
	$$\mathbb{P}(||Z|| > \varpi \sqrt{n}) \leq c^* \exp(-c \varpi n)$$ 
	for all  $\varpi \in \mathbb{R}$ greater than or equal to $c^*$.
\end{lem}
As a consequence, we have $||Z|| \leq \varpi \sqrt{n}$ with high probability.
\begin{definition}[Definition 5.1 of \cite{vershyninintroduction}, Nets]
	Let $(X, d)$ be a metric space and $\delta > 0$. A subset $\mathcal{N}_\delta$ of $X$ is called a $\delta$-net of $X$ if every point $x \in X$ can be approximated to within $\delta$ by some point $y \in \mathcal{N}_\delta$, such that $d(x, y) \leq \delta$.
\end{definition}
\begin{lem}[Lemma 5.2 of \cite{vershyninintroduction}]\label{vers5.2}
	The unit Euclidean sphere\footnote{The unit Euclidean $(n-1)$-sphere is defined as $\mathbb{S}^{n - 1} = \{x \in \mathbb{R}^n : ||x|| = 1$\}} $\mathbb{S}^{n - 1}$ equipped with the Euclidean metric satisfies, for every $\delta > 0$, $\mathcal{N}_ \delta \leq \left(1 + \frac{2}{\delta}\right)^n.$
\end{lem}

\begin{lem}[Lemma 5.3 of \cite{vershyninintroduction}]\label{vers5.3}
	Let $Z \in \mathbb{R}^{n \times n}$ be a matrix, and suppose that $\mathcal{N}_\delta$ is a $\delta - net$ of $\mathbb{S}^{n - 1}$ for some $\delta \in [0, 1)$. Then
	$$||Z|| \leq (1 - \delta)^{-2} \sup_{x, y \in \mathcal{N}_\delta} \langle y, Z x \rangle .$$
\end{lem}
Equipped with these Lemmas, we are now ready to prove Lemma \ref{mymainlemma}. 
\subsection{Proof of Lemma \ref{mymainlemma}}
Proof of (i.)
\begin{proof}
	Note that from the construction of the dual certificate, $Q = Q_B + Q_C$. Also, from Equation (\ref{golfingscheme}), $Y_k = Y_{k - 1} + p^{-1}\mathcal{P}_{\Gamma_k} \mathcal{P}_T(UU^T - Y _{k - 1})$ and $Y_0 = 0$. We set 
	\begin{equation}\label{zkrecursion}
	Z_k = UU^T - \mathcal{P}_T Y_k \text{ for } k = 0, \ldots, k_0;
	\end{equation}
	 so that $Z_0 = UU^T$ and $Z_k = (\mathcal{P}_T -  p^{-1}\mathcal{P}_T \mathcal{P}_{\Gamma _k} \mathcal{P}_T)Z_{k-1}.$
	 
	Hence, 
	\begin{equation}\label{ykzk}
	Y_{k_0} = \sum_{k = 1}^{k_0}p^{-1} \mathcal{P}_{\Gamma_k}Z_{k-1}.
	\end{equation}
	Observe that $Z_k \in T$, therefore $\mathcal{P}_{T^\perp}Z_k = 0$. Since $Q_B = \mathcal{P}_{T^\perp}Y_{k_0}$ (see Equation \eqref{eq:qbequalyk0}), then
	\begin{align}
	||Q_B|| = ||\mathcal{P}_{T^\perp}Y_{k_0}|| &\leq \sum_{k = 1}^{k_0}||p^{-1}\mathcal{P}_{T^\perp} \mathcal{P}_{\Gamma_k} Z_{k - 1}|| \nonumber \\ 
	&= \sum_{k = 1}^{k_0}||\mathcal{P}_{T^\perp}(p^{-1} \mathcal{P}_{\Gamma_k} Z_{k - 1} - Z_{k - 1})||^\diamondsuit \nonumber \\
	& \leq \sum_{k = 1}^{k_0}||p^{-1} \mathcal{P}_{\Gamma_k} Z_{k - 1} - Z_{k - 1}|| \nonumber \\
	& = \sum_{k = 1}^{k_0}||(p^{-1} \mathcal{P}_{\Gamma_k} - \mathcal{I})Z_{k - 1}|| \nonumber \\
	& \leq c\sum_{k = 1}^{k_0}\left(\frac{\log n}{p}||Z_{k - 1}||_\infty +  \sqrt{\frac{\log n}{p}}||Z_{k - 1}||_{\infty, 2} \right) ^\spadesuit  \nonumber \\
	& \leq c \sum_{k = 1}^{k_0} \left(\frac{n}{c_0 \mu r}||Z_{k - 1}||_\infty + \sqrt{\frac{n}{c_0 \mu r}}||Z_{k - 1}||_{\infty, 2} \right) ^ \clubsuit \nonumber \\
	& \leq c \sum_{k = 1}^{k_0} \left(\frac{n}{\sqrt{c_0 \mu r}}||Z_{k - 1}||_\infty + \sqrt{\frac{n}{c_0 \mu r}}||Z_{k - 1}||_{\infty, 2} \right) ^ \dagger \nonumber \\
	& \leq \frac{c}{\sqrt{c_0}} \sum_{k = 1}^{k_0} \left(\frac{n}{\sqrt{\mu r}}||Z_{k - 1}||_\infty + \sqrt{\frac{n}{\mu r}}||Z_{k - 1}||_{\infty, 2} \right). \label{dagadaga}
	\end{align} 
	The expression in the first line is as a result of Equation \eqref{ykzk}. $\diamondsuit$ follows from the fact that $\mathcal{P}_{T^\perp} Z_k = 0$. $\spadesuit$ is application of Lemma \ref{lemma2yudongchen}. $\clubsuit$  holds for $p \geq \frac{c_0 \mu  r \log n}{n}$ in Lemma \ref{lemma4chenincoherence} while $\dagger$ is valid since we can choose $c_0$ such that $c_0 \mu r > 1$ so that $\frac{n}{c_0 \mu r} \leq \frac{n}{\sqrt{{c_0 \mu r}}}$ holds. We will now bound $||Z_{k - 1}||_\infty$ and $||Z_{k - 1}||_{\infty, 2}$ by applying Lemma \ref{lemma4chenincoherence}. Using \eqref{zkrecursion}, applying Lemma \ref{lemma4chenincoherence} repeatedly (replacing $\Gamma$ with $\Gamma_k$), we have (with high probability);
	\begin{align}
	||Z_{k - 1}||_\infty &= ||(\mathcal{P}_T - p^{-1}\mathcal{P}_T \mathcal{P}_{\Gamma_{k - 1}}\mathcal{P}_T) \ldots (\mathcal{P}_T - p^{-1}\mathcal{P}_T \mathcal{P}_{\Gamma_1}\mathcal{P}_T)Z_0||_\infty \nonumber \\
	&= \left|\left| \prod _{k = 1} ^{k - 1} (\mathcal{P}_T - p^{-1}\mathcal{P}_T \mathcal{P}_{\Gamma_{k}}\mathcal{P}_T) Z_0\right|\right|_\infty \nonumber \\
	&\leq (1/2)^{k - 1} ||UU^T||_ \infty,\label{star} 
	\end{align}
	since $Z_0 = UU^T$. In like manner, we apply Lemma \ref{lemma3chencoherence} to $||Z_{k - 1}||_{\infty, 2}$, using Equation \eqref{zkrecursion} again and replacing $\Gamma$ with $\Gamma_k$ to get, with high probability,
	\begin{align}
	||Z_{k - 1}||_{\infty, 2} &= ||(\mathcal{P} _ T - p^{-1}\mathcal{P} _T \mathcal{P} _ {\Gamma_{k - 1}} \mathcal{P}_T)Z_{k - 2}||_{\infty, 2} \nonumber \\
	&\leq 1/2 \sqrt{\frac{n}{\mu r}}||Z_{k - 2}||_\infty + 1/2 ||Z_{k - 2}||_{\infty, 2} \label{starstar}
	\end{align}
	Combining (\ref{star}) and (\ref{starstar}) and using (\ref{zkrecursion}) repeatedly, we have
	\begin{align}
	&||Z_{k - 1}||_{\infty, 2} \leq 1/2 \sqrt{\frac{n}{\mu r}} \Vert Z_{k - 2} \Vert _\infty + 1/2 \Vert Z_{k - 2} \Vert _{\infty , 2} \nonumber \\ \nonumber
	&\leq 1/2 \sqrt{\frac{n}{\mu r}} \Vert Z_{k - 2} \Vert _\infty + 1/2 \left(1/2 \sqrt{\frac{n}{\mu r}} \Vert Z_{k - 3} \Vert _\infty  + 1/2 \Vert Z_{k - 3} \Vert _{\infty , 2}\right)\\\nonumber
&\leq  \sqrt{\frac{n}{\mu r}} \left( 1/2 \Vert Z_{k - 2} \Vert _\infty + 1/4 \Vert Z_{k - 3} \Vert _\infty  + 1/8  \Vert Z_{k - 4} \Vert \right) + 1/8 \Vert Z_{k - 4} \Vert _{\infty , 2}\\	\nonumber
&\leq \sqrt{\frac{n}{\mu r}}\left(\frac{1}{2} \times \frac{1}{2^{k-2}} + \frac{1}{4} \times \frac{1}{2^{k-3}} + \cdots + \frac{1}{2^{k-1}}\right)\Vert UU^T \Vert _ \infty +  \frac{1}{2^{k-1}} \Vert UU^T \Vert _ {\infty,2}\\
	&\leq k(1/2)^{k - 1} \sqrt{\frac{n}{\mu r}}||UU^T||_\infty + (1/2)^{k - 1}||UU^T||_{\infty, 2}. \label{zinf2}
	\end{align}
	Substituting (\ref{star}) and (\ref{zinf2}) back into (\ref{dagadaga}), we get;
	\begin{align*}
	||Q_B|| &\leq \frac{c}{\sqrt{c_0}} \sum_{k = 1}^{k_0} \left(\frac{n}{\sqrt{\mu r}}(1/2)^{k - 1}||UU^T||_\infty + \right. \\
	&\left. \sqrt{\frac{n}{\mu r}}\left(k(1/2)^{k - 1}\sqrt{\frac{n}{\mu r}} ||UU^T||_\infty + (1/2)^{k - 1}||UU^T||_{\infty, 2} \right) \right)\\
	&= \frac{c}{\sqrt{c_0}} \sum_{k = 1}^{k_0} \left(\frac{n}{\sqrt{\mu r}}(k + 1)(1/2)^{k - 1}||UU^T||_\infty + \sqrt{\frac{n}{\mu r}}(1/2)^{k - 1}||UU^T||_{\infty, 2} \right)\\
	&\leq \frac{c}{\sqrt{c_0}} \frac{n}{\sqrt{\mu r}} ||UU^T||_\infty \sum_{k = 1}^{k_0}(k + 1)(1/2)^{k - 1} + \frac{c}{\sqrt{c_0}} \sqrt{\frac{n}{\mu r}} ||UU^T||_{\infty, 2} \sum_{k = 1}^{k_0}(1/2)^{k - 1}\\
	&\leq \frac{6c}{\sqrt{c_0}} \frac{n}{\sqrt{\mu r}} ||UU^T||_\infty + \frac{2c}{\sqrt{c_0}} \sqrt{\frac{n}{\mu r}} ||UU^T||_{\infty, 2}.
	\end{align*}
	The constants $6$ and $2$ on the last line above are the partial sums of the sequences in the previous line. From (\ref{jointincoherence}) and (\ref{incoherenceforlinfinity2norm}), $||UU^T||_\infty \leq \frac{\sqrt{\mu r}}{n}$ and 
	$||UU^T||_{\infty, 2}  \leq \sqrt{\frac{\mu r}{n}}.$
	So, we have
	\begin{align*}
	||Q_B|| = ||\mathcal{P}_{T^\perp} Y_{k_0}|| & \leq  \frac{6c}{\sqrt{c_0}} + \frac{2c}{\sqrt{c_0}}\\
	&= \frac{8c}{\sqrt{c_0}} \leq 1/8,
	\end{align*}
	provided $c_0 \geq (64c)^2$.
\end{proof}

Proof of (ii.)
\begin{proof}
	\begin{align}
	\mathcal{P}_\Gamma (UU^T + Q_B) &= \mathcal{P}_\Gamma (UU^T + \mathcal{P}_{T^\perp} Y_{k_0}) \nonumber \\
	&= \mathcal{P}_\Gamma (UU^T + Y_{k_0} - \mathcal{P}_{T} Y_{k_0})\nonumber \\
	&= \mathcal{P}_\Gamma (UU^T - \mathcal{P}_{T} Y_{k_0}) \label{bigstar} \\
	&= \mathcal{P}_\Gamma (UU^T - \mathcal{P}_{T}(Y_{k_0 - 1} + p^{-1}\mathcal{P}_{\Gamma{k_0}} \mathcal{P}_T(UU^T -  Y_{k_0 - 1}))) \nonumber \\
	&= \mathcal{P}_\Gamma (UU^T - \mathcal{P}_{T} Y_{k_0 - 1} -  p^{-1}\mathcal{P}_T \mathcal{P}_{\Gamma_{k_0}} \mathcal{P}_T(UU^T -  Y_{k_0 - 1})) \nonumber \\
	&= \mathcal{P}_\Gamma (Z_{k_0 - 1} - p^{-1}\mathcal{P}_T \mathcal{P}_{\Gamma_{k_0}} \mathcal{P}_T Z_{k_0 - 1}) \nonumber \\
	&= \mathcal{P}_\Gamma (\mathcal{P}_T -p^{-1} \mathcal{P}_T \mathcal{P}_{\Gamma_{k_0}} \mathcal{P}_T) Z_{k_0 - 1} = \mathcal{P}_\Gamma (Z_{k_0}). \nonumber
	\end{align}
	$\mathcal{P}_\Gamma Y_{k_0} = 0$ (since $\Gamma \cap T = \{0\}$ from the conditions of Lemma \ref{lemma2.4candes}), the fourth line is as a result of Equation \eqref{golfingscheme}, the sixth line follows from the definition of $Z_k$ in Equation \eqref{zkrecursion} above;  therefore \eqref{bigstar} holds. Hence,
	\begin{equation}\label{proofzk01}
	||\mathcal{P}_\Gamma (UU^T + Q_B)|| \leq ||Z_{k_0}|| \leq ||\mathcal{P}_T - p^{-1}\mathcal{P}_T \mathcal{P}_{\Gamma_{k_0}} \mathcal{P}_T)||.||Z_{k_0 - 1}||.
	\end{equation}
	Notice that by the conditions of Theorem \ref{main}, $p \geq \frac{c_0 \mu r \log n}{n}$ and $\Gamma_k$ is independent of $Z_{k - 1}$. Using the recurrence relation \eqref{zkrecursion} again and applying Lemma \ref{lemma1chenchoerence} repeatedly (replacing $\Gamma$ with $\Gamma_k$), we get with high probability;
	\begin{equation}\label{proofzk02}
	\begin{aligned}
	||\mathcal{P}_\Gamma (UU^T + Q_B)||_\infty &\leq (1/2)^{k_0} ||UU^T||_\infty \\
	&= (1/2)^{k_0}  \sqrt{\frac{\mu r}{n^2}} \\
	& < \frac{1}{8\sqrt{n}} = \frac{\lambda}{8}, (\text{we chose } \lambda = \frac{1}{\sqrt{n}}).
	\end{aligned}
	\end{equation}
The last inequality is valid for $n \geq 16^2 \mu r$ since the term $(1/2)^{k_0}$ is small for large $k_0$.
\end{proof}

Proof of (iii.) 
\begin{proof} 
	We are required to show that $||\mathcal{P}_{\Gamma^C} (UU^T + Q_B)||_\infty < \lambda / 4$. Observe that $Z_k = UU^T - \mathcal{P}_T Y_k$ implies that $UU^T = Z_k + \mathcal{P}_T Y_k$. Hence,
	\begin{align*}
	UU^T + Q_B &= Z_{k_0} + \mathcal{P}_T Y_{k_0} + \mathcal{P}_{T^\perp}Y_{k_0}\\
	&= Z_{k_0} + Y_{k_0}.
	\end{align*}
	We know that the support of $Y_{k_0}$ is $\Gamma^C$ and, by virtue of \eqref{proofzk01} and \eqref{proofzk02}, we have shown that $||Z_{k_0}|| = ||\mathcal{P}_T - p^{-1} \mathcal{P}_T \mathcal{P}_{\Gamma_{k_0}} \mathcal{P}_T) Z_{k_0 - 1}|| < \lambda /8$. Therefore, to show that $||\mathcal{P}_{\Gamma^C}(UU^T + Q_B)||_\infty = ||Z_{k_0} + Y_{k_0}||_\infty < \lambda / 4$, it is sufficient to show that $||Y_{k_0}||_\infty < \lambda / 8$. Recall from \eqref{ykzk} that $Y_{k_0} = \sum_{k = 1}^{k_0} p^{-1}\mathcal{P}_{\Gamma_k}Z_{k - 1}$. It then follows that
	\begin{align*}
	||Y_{k_0}||_\infty &\leq  p^{-1}\sum_{k = 1}^{k_0} ||\mathcal{P}_{\Gamma_k}Z_{k - 1}||_\infty\\
	&\leq p^{-1} \sum_{k = 1}^{k_0} ||Z_{k - 1}||_\infty\\
	& \leq p^{-1} \sum_{k = 1}^{k_0}(1/2)^{k - 1}||UU^T||_\infty\\
	&\leq p^{-1} \sum_{k = 1}^{k_0}(1/2)^{k - 1} \frac{\sqrt{\mu r}}{n}\\
	&\leq \frac{n}{c_0 \mu r \log n} \sum_{k = 1}^{k_0}(1/2)^{k - 1} \frac{\sqrt{\mu r}}{n}\\
	&\leq \frac{1}{c_0 \sqrt{\mu r} \log n} \sum_{k = 1}^{k_0}(1/2)^{k - 1} < \frac{1}{8\sqrt{n}} = \frac{\lambda}{8}.
	\end{align*}
	The third line above is application of Equation \eqref{star} to the sequence, the fourth line is application of the coherence property \eqref{jointincoherence}, the first inquality on the last line follows from the conditions of Theorem \ref{main} while the  second inequality holds for $n \geq \exp(16/c_0\sqrt{r})$.
\end{proof}
Proof of (iv.) 
\begin{proof}
	From equation (\ref{defq_c_alternative}),
	\begin{align*}
	Q_C &= \lambda \mathcal{P} _{T^\perp} \sum_{k \geq 0} (\mathcal{P}_\Gamma \mathcal{P} _T \mathcal{P} _\Gamma)^k Sgn(C)\\
	&= \lambda \mathcal{P}_{T^\perp} Sgn(C) + \lambda \mathcal{P}_{T^\perp} \sum_{k \geq 1} (\mathcal{P}_\Gamma \mathcal{P} _T \mathcal{P} _\Gamma)^k Sgn(C).
	\end{align*}
	We define $\Phi = Sgn(C)$ with the following probability distribution:
	\begin{equation*}
	\Phi_{ij} = \begin{cases}
	1 \qquad \text{with probability } p / 2, \\
	0 \qquad \text{with probability } 1 - p, \\
	-1 \quad \text{ with probability } p / 2. \\
	\end{cases}
	\end{equation*}
	Then,
	\begin{align}
	Q_C &= \lambda \mathcal{P}_{T^\perp} \Phi + \lambda \mathcal{P}_{T^\perp} \sum_{k \geq 1} (\mathcal{P}_\Gamma \mathcal{P} _T \mathcal{P} _\Gamma)^k \Phi \nonumber \\
	& \coloneqq \mathcal{P} _{T^\perp} Q_C^0 + \mathcal{P} _{T^\perp} Q_C^1. \label{plusstar}
	\end{align}
	Therefore, to derive a bound for $Q_C$, we only need to derive a bound for  $\mathcal{P} _{T^\perp} Q_C^0$ and $\mathcal{P} _{T^\perp} Q_C^1$.  Applying Lemma \ref{terry235} to the first term of (\ref{plusstar}), the following holds with high probability;
	\begin{align*}
	|| \mathcal{P} _{T^\perp} Q_C^0|| & \leq |\lambda|||\Phi|| \\
	& \leq \lambda c^* \sqrt{p n} = c^* \sqrt{p}.
	\end{align*}
	The bound on the spectral norm of $\Phi$ is as a result of Lemma \ref{terry235} and recall that $\lambda = \frac{1}{\sqrt{n}}$, hence the last part of the equation above follows. For the second term of (\ref{plusstar}), let $\mathcal{R} = \mathcal{P} _{T^\perp} \sum_{k \geq 1} (\mathcal{P}_\Gamma \mathcal{P}_T \mathcal{P}_\Gamma)^k$, so that $||Q_C^1|| = |\lambda|||\mathcal{R}(\Phi)||$. We are required to find a bound on the operator norm of $\mathcal{R}(\Phi)$. We make use of Lemmas \ref{vers5.2} and \ref{vers5.3} to achieve this. Let $\mathcal{N}_\delta$ be a $\delta - net$ for a sphere $\mathbb{S}^{n - 1}$ of size at most $(3/ \delta)^n$. Then by Lemma \ref{vers5.3},
	\begin{equation*}
	||\mathcal{R}(\Phi)|| \leq (1 - \delta)^{-2} \max_{x , y \in \mathcal{N}_\delta} \langle y, \mathcal{R}(\Phi)x \rangle
	\end{equation*}
	We define a random variable $Y(x,y)$ such that
	\begin{align}
	Y(x,y) & \coloneqq \langle y, \mathcal{R} (\Phi) x \rangle \nonumber \\
	&= \langle \mathcal{R}(yx^*), \Phi \rangle \label{dagger}
	\end{align}
	for a fixed pair $(x, y) \in \mathcal{N}_\delta \times \mathcal{N}_\delta$ with unit norm. The equality in Equation \eqref{dagger} holds because $\mathcal{R}$ is a self-adjoint operator \cite{candesrobust}. Note that the support of $\Phi$ is $\Gamma$ and $\Phi$ is a symmetric matrix with independent and identically distributed random signs. Using Hoeffding's Inequality, conditioning on $\Gamma$, we have
	\begin{equation*}
	\mathbb{P}(|Y(x,y)| > \tau | \Gamma) \leq 2 \exp \left(\frac{-2 \tau^2}{||\mathcal{R}(yx^*)||^2_F}\right),
	\end{equation*}
	for $\tau \geq 0$.
	Since $(yx^*)$ is a unit-normed vector, $||yx^*||_F = 1$. Therefore, $||\mathcal{R}(yx^*)||_F \leq ||\mathcal{R}||$ holds. So, we have
	\begin{equation}\label{hoefding_bound_for_R}
	\mathbb{P}\left( \max_{x,y \in \mathcal{N}_\delta} |Y(x,y)| > \tau | \Gamma \right) \leq 2 |\mathcal{N}_\delta|^2 \exp \left( \frac{-2 \tau^2}{||\mathcal{R}||^2}\right).
	\end{equation}
	\begin{align*}
	||\mathcal{R}|| &= ||\mathcal{P}_{T^\perp} \sum_{k \geq 1} (\mathcal{P}_\Gamma \mathcal{P}_T \mathcal{P}_\Gamma)^k||\\
	& \leq \sum_{k \geq 1} ||(\mathcal{P}_\Gamma \mathcal{P}_T \mathcal{P}_\Gamma)^k ||  = \sum_{k \geq 1} ||(\mathcal{P}_\Gamma \mathcal{P}_T)||^{2k} \\
	&= \frac{||\mathcal{P}_\Gamma \mathcal{P}_T||^2}{1 - ||\mathcal{P}_\Gamma \mathcal{P}_T||^2}.
	\end{align*}
	The equality in the immediate equation before the last line above holds with the assumption that the product of $\mathcal{P} _ T$ and $\mathcal{P} _ \Gamma$ commute and for the fact that $\mathcal{P} _ \Gamma ^ 2 = \mathcal{P} _ \Gamma$. Considering the event $\psi \coloneqq \{||\mathcal{P}_\Gamma \mathcal{P}_T|| \mid ||\mathcal{P}_\Gamma \mathcal{P}_T|| \leq \beta\}$,
	\begin{equation*}
	||\mathcal{R}|| \leq \frac{||\mathcal{P}_\Gamma \mathcal{P}_T||^2}{1 - ||\mathcal{P}_\Gamma \mathcal{P}_T||^2} = \frac{\beta^2}{1 - \beta^2}
	\end{equation*}
	Now, replacing $\max_{x,y \in \mathcal{N}_\delta} |Y(x,y)|$ in Equation \eqref{hoefding_bound_for_R} with $||\mathcal{R}(\Phi)||$ and substituting $|\mathcal{N}_\delta|$ from Lemma \eqref{vers5.2} together with almost sure bound of $\mathcal{R}(\Phi)$, we get
	\begin{align*}
	\mathbb{P}\left(||\mathcal{R}(\Phi)|| > \tau |\Gamma \right) &\leq 2 ( 3/\delta)^{2n} \exp \left( \frac{-2 \tau^2}{((1 - \delta)^{-2}(\frac{\beta^2}{1 - \beta^2}))^2}\right)\\
	&= 2 (3/\delta)^{2n} \exp\left(\frac{-2\tau^2 (1-\beta^2)^2(1 - \delta)^4}{\beta^4}\right),
	\end{align*}
	and marginally,
	\begin{equation*}
	\mathbb{P}\left(||\mathcal{R}(\Phi)|| > \tau \right) \leq 2 (3/\delta)^{2n} \exp\left(\frac{-2\tau^2 (1-\beta^2)^2(1 - \delta)^4}{\beta^4}\right) +  \mathbb{P}(||\mathcal{P}_\Gamma \mathcal{P} _ T|| \geq \beta).
	\end{equation*}
	Consequently, 
	\begin{equation*}
	\mathbb{P}\left(\lambda ||\mathcal{R}(\Phi)|| > \tau \right) \leq 2 (3/\delta)^{2n} \exp\left(\frac{-2\tau^2 (1-\beta^2)^2(1 - \delta)^4}{\lambda^2 \beta^4}\right) +  \mathbb{P}(||\mathcal{P}_\Gamma \mathcal{P} _ T|| \geq \beta),
	\end{equation*}
	with $\lambda = \frac{1}{\sqrt{n}}$.
	Conclusively,
	$$||Q_C|| \leq   c^* \sqrt{p} + |\lambda|||\mathcal{R}(\Phi)||  \leq 1/8,$$
	with high probability, provided that $\beta$ is small enought and for an appropriate choice of $c^*$ and  $\delta$. 
\end{proof}
Proof of (v.)
\begin{proof}
	Recall that $Q_C = \lambda \mathcal{P} _{T^\bot}(\mathcal{P} _ \Gamma - \mathcal{P} _ \Gamma \mathcal{P} _ T \mathcal{P} _ \Gamma)^{-1}Sgn(C)$. Therefore,
	\begin{align*}
	\mathcal{P} _ {\Gamma^\bot} Q_C &=  \lambda \mathcal{P} _ {\Gamma^\bot} \mathcal{P} _{T^\bot}(\mathcal{P} _ \Gamma - \mathcal{P} _ \Gamma \mathcal{P} _ T \mathcal{P} _ \Gamma)^{-1}Sgn(C)\\
	&=  \lambda \mathcal{P} _ {\Gamma^\bot}(\mathcal{I} -  \mathcal{P} _{T})(\mathcal{P} _ \Gamma - \mathcal{P} _ \Gamma \mathcal{P} _ T \mathcal{P} _ \Gamma)^{-1}Sgn(C)\\
	&=  \lambda [\mathcal{P} _ {\Gamma^\bot} (\mathcal{P} _{T}\mathcal{P} _ \Gamma - \mathcal{P} _ \Gamma \mathcal{P} _ T \mathcal{P} _ \Gamma)^{-1} -  \mathcal{P} _{T}(\mathcal{P} _ \Gamma - \mathcal{P} _ \Gamma \mathcal{P} _ T \mathcal{P} _ \Gamma)^{-1}]Sgn(C).
	\end{align*}
	Since the operator $\mathcal{P} _ \Gamma - \mathcal{P} _ \Gamma \mathcal{P} _ T \mathcal{P} _ \Gamma$ maps $\Gamma$ onto itself, $(\mathcal{P} _ \Gamma - \mathcal{P} _ \Gamma \mathcal{P} _ T \mathcal{P} _ \Gamma)^{-1}$ is in $\Gamma$ and $\mathcal{P} _ {\Gamma^\bot}(\mathcal{P} _ \Gamma - \mathcal{P} _ \Gamma \mathcal{P} _ T \mathcal{P} _ \Gamma)^{-1} = 0$. So, we have
	\begin{equation}
	\mathcal{P} _ {\Gamma^\bot} Q_C = -  \lambda  \mathcal{P} _ {\Gamma^\bot} \mathcal{P} _{T} (\mathcal{P} _ \Gamma - \mathcal{P} _ \Gamma \mathcal{P} _ T \mathcal{P} _ \Gamma)^{-1}\Phi,
	\end{equation}
	where $\Phi = Sgn(C)$. Now, for $(i,j) \in \Gamma^C$,
	\begin{align*}
	{Q_C} _{ij} &= \langle e_i, Q_C e_j \rangle\\
	&= \langle e_i e_j^T, Q_C \rangle\\
	&= \lambda \langle W(i,j), \Phi \rangle,
	\end{align*}
	where $W(i,j)$ is the matrix $-(\mathcal{P} _ \Gamma - \mathcal{P} _ \Gamma \mathcal{P} _ T \mathcal{P} _ \Gamma)^{-1} \mathcal{P} _ \Gamma \mathcal{P} _ T (e_i e_j^T)$. Since $\Gamma = supp(\Phi)$, the signs of $\Phi$ are independent, symmetric and identically distributed. Applying Hoeffding's inequality, we have
	$$\mathbb{P}(|{Q_C} _{ij}| > \tau \lambda | \Gamma) \leq 2 exp \left( \frac{-2 \tau^2}{||W(i,j)|| _F ^2}\right),$$
	therefore,
	$$\mathbb{P}(\sup_{i,j} |{Q_C} _{ij}| > \tau \lambda | \Gamma) \leq 2 n^2 exp \left( \frac{-2 \tau^2}{\sup _{i,j} ||W(i,j)|| _F ^2}\right).$$
	For a matrix of the form $e_ie_j^T$, from \eqref{equ:mat_ort_pro}, we have
	\begin{equation}
	||\mathcal{P} _{T^\bot} e_i e_j^T||_F^2 = ||(\mathcal{I} - UU^T)e_i||^2||(\mathcal{I} - UU^T)e_j||^2 \geq \left(1 - \frac{\mu r}{n}\right)^2.
	\end{equation}
	The last inequality is based on \eqref{standardincoherence} with the assumption that $\frac{\mu r}{n} \leq 1$. From the fact that $||\mathcal{P} _ T e_i e_j^T||_F^2 + ||\mathcal{P} _{T^\bot} e_i e_j^T||_F^2 = 1$, we deduce 
	\begin{equation}
	||\mathcal{P} _ T e_i e_j^T||_F \leq \sqrt{\frac{2 \mu r}{n}},
	\end{equation}
	so that
	\begin{align*}
	||\mathcal{P}_\Gamma \mathcal{P} _ T (e_i e_j^T)||_F & \leq ||\mathcal{P} _ \Gamma \mathcal{P} _ T||||\mathcal{P}_T(e_ie_j^T)||_F\\
	&\leq \beta \sqrt{\frac{2 \mu r}{n}},
	\end{align*}
	on the event $\{||\mathcal{P}_\Gamma \mathcal{P}_T|| \leq \beta \}$. Similarly, on the same event,
	\begin{equation}
	||(\mathcal{P} _ \Gamma - \mathcal{P} _ \Gamma \mathcal{P} _ T \mathcal{P} _ \Gamma)^{-1}||_F \leq \frac{1}{1 - \beta ^ 2}.
	\end{equation}
	Therefore,
	\begin{equation}
	||W(i,j)||^2_F \leq \frac{2 \beta ^2}{(1 - \beta^2)^2} \frac{\mu r}{n}.
	\end{equation}
	Hence, unconditionally, we have
	\begin{equation}
	\mathbb{P}\left(\sup_{i,j} |{Q_C} _{ij}| > \tau \lambda\right) \leq 2 n^2 exp \left( \frac{- n (1 - \beta^2)^2 \tau^2}{2 \beta^2 \mu r}\right) + \mathbb{P}(||\mathcal{P}_\Gamma \mathcal{P} _ T|| \geq \beta).
	\end{equation}
\end{proof}
Similar proof to this part can be found in \cite{candesrobust}. This concludes the proof of the Theorems. We performed series of numerical experiments to corroborate our claim. The report of these experiments is presented in the next section.

\section{Experiments}\label{experiments}
In this section, we report the results of the numerical experiments performed to test the efficiency and effectiveness of problem \eqref{qcp1} in planted maximum quasi-clique recovery. Firstly, we investigated the performance of the problem for various $\gamma$-clique and graph sizes. Secondly, we fixed the size of the graph and the planted $\gamma$-clique and then varied the edge connectivity within the $\gamma$-clique and the amount of diversionary edges that can be added for which perfect recovery would still be realized. All of these experiments were performed using CVXPY developed by \citeauthor{cvxpy} \cite{cvxpy}. The experiments were performed in Python on an HP personal computer equipped with Intel core $i7$ 2.7 GHz processor and 16GB RAM. For both of the experiments, we chose $\lambda = \frac{1}{\sqrt{n}}$. There is nothing special about our choice of $\lambda$. Various other values of $\lambda$ also works \cite{candesrobust,chenincoherence}. Each of the experiments was repeated ten times.

\subsection{Exact recovery from varying $\gamma-clique$ size}\label{sec:72}
For our model to recover $\gamma$-clique, $n_c$ must be large enough. To check what size of $\gamma$-clique can be recovered from a given network, we planted $\gamma$-clique of various sizes in some graphs and tried to recover them. We randomly considered $n \times n$  symmetric binary matrices  with $n = 25, 50,75, \ldots 200$. This gives eight different matrix sizes. For each of these different matrices, we inserted $\gamma$-clique of sizes equal to $10\%, 20\%, \dots , 100\%$ of the size of the matrices using binomial distribution as described before. In total, we have eighty problem instances and each instance was repeated ten times, giving rise to $800$ runs in total. In each case, we fixed $\gamma$ (the edge density parameter of the quasi-clique) to be $0.85$ and $\rho$ (the probability of adding a diversionary edge) was $0.25$. The details of the matrix construction is as follow. We generated an $n \times n$ zero matrices $B \in \mathbb{R} ^{n \times n}$ (the adjacency matrix of the planted quasi clique) and $C \in \mathbb{R} ^{n \times n}$ (the adjacency matrix of the nodes not belonging to the quasi clique). Let $\Gamma$ be the set of indices corresponding to the nodes belonging to the quasi-clique. Consequently, $\Gamma ^c$ is the set of indices of the nodes of the diversionary edges. The entries of $\Gamma$ are equal to one with probability $\gamma$ and zero with probability $1 - \gamma$. Hence, $\Gamma$ has expected cardinality $|\Gamma| = \gamma n_c^2$, where $n_c$ is the dimension of the non-zero submatrix of $B$ corresponding to the adjacency matrix of the planted quasi clique. $\Gamma ^c$ also follows a Bernoulli distribution with parameter $\rho$ and expected cardinality of $\rho (n^2 - n_c^2)$. 
We formed $A$ by adding $B$ and $C$ and tried to recover $B$ from $A$. We declare a recovery attempt to be successful if the relative error (in Frobenius norm) of the recovered matrix with respect to $B$ is less than or equal to $10^{-6}$, i.e, $||B - B^*||_F/||B||_F \leq 10^{-6}$, where $B^*$ is the recovered $\gamma$-clique. 

\begin{figure}[]
	\centering
	\includegraphics[width=0.75\textwidth]{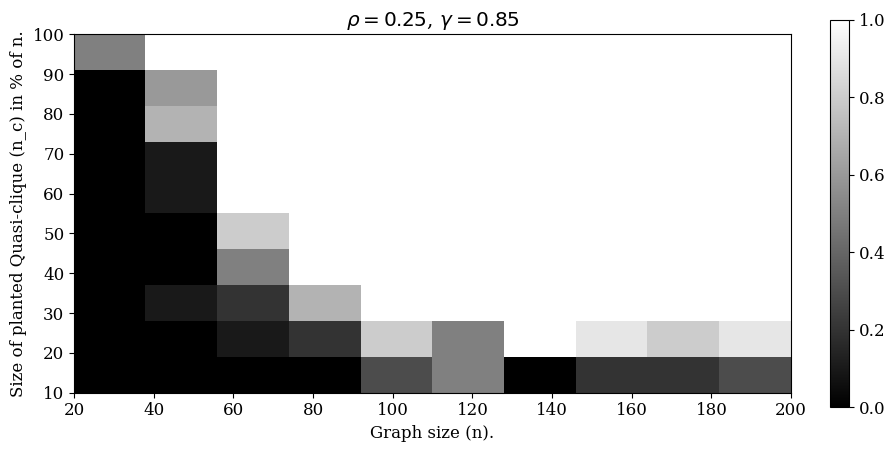}
	\caption{Exact recovery of varying quasi-clique size from graphs of different sizes \label{qvsn}}
\end{figure}

Figure \ref{qvsn} contains the plot of the probability of recovery from different graph sizes and quasi-clique sizes. White area denotes perfect recovery while the black area means the recovery attempt failed for every trial. We discovered that when the $\gamma$-clique is less than $10\%$ of the graph size, recovery is impossible for all the graph sizes  considered. In addition, when the quasi-clique size is bigger than $60\%$ of the graph size, perfect recovery is guaranteed for all the graph sizes in our experiment. We also observed that as graph size ($n$) increases, the minimum proportion of the size of planted quasi clique $(n_c)$, with respect to $n$, required for perfect recovery reduces.

\subsection{Recovery from varying edge density and random noise}\label{sec:73}
\begin{figure}[]
	\begin{center}
		\includegraphics[width =0.75\textwidth]{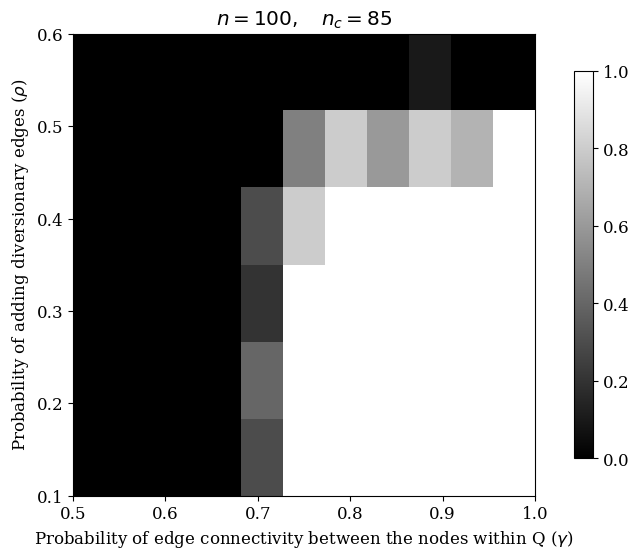}
		\caption{$\gamma - clique$ recovery from varying edge density and random noise with $n = 100$ and $n_c = 85$ \label{rhovsgamma100}}
	\end{center}
\end{figure}

\begin{figure}[H]
	\begin{center}
		\includegraphics[width =0.75\textwidth]{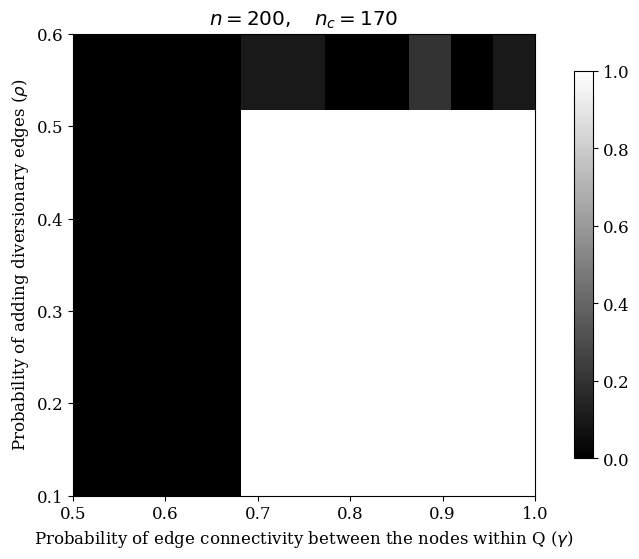}
		\caption{$\gamma$-clique recovery from varying edge density and random noise with $n = 200$ and $n_c = 170$ \label{rhovsgamma200}}
	\end{center}
\end{figure}

Theorem \ref{main} states that convex programming can recover an incoherent, low-rank noisy matrix with a reasonable level of missing entries. In the next experiment, we determined, experimentally, the maximum tolerable level of diversionary edges (noise) that can be added along with the minimum (edge density) probability within the planted quasi clique that is required for a perfect recovery. The set-up follows the previous experiment with $n = 100 \text{ and } 200$. In both cases, the $\gamma$-clique size, $n_c$, was $85$ and $170$ respectively. The results of this experiment is presented in Figures \ref{rhovsgamma100} and \ref{rhovsgamma200}. In both figures, the probability of the edge connectivity within the planted quasi clique ($\gamma$) is plotted on the $x$-axis against the probablity of adding a diversionary edge on the $y$-axis. As above, the heatmap shows the probablity of recovery for each instance of the experiment. White area shows perfect recovery while black implies no recovery. In both cases, we observed a phase transition, whereby the probability of recovery jumps from zero to one as $\gamma$ reached a particular threshold. This is in agreement with our theoretical claim that when $\gamma \geq \frac{c_0 \mu r \log n}{n}$, for some constant $c_0 >0$, exact recovery is guaranteed. In addition, we recorded complete failure when the probability of adding a random noise, $\rho,$ is greater than or equal to $0.6$, i.e, when $\rho \geq 0.6$. This implies that when the random noise becomes too much, recovery becomes impossible.  As a conclusion, no $\gamma$-clique of any size can be recovered by our method when $\gamma$ is less than the recovery threshold (i.e, the edge density of the planted $\gamma$-clique to be recovered is too low), the size of the planted $\gamma$-clique is too small or the probability of adding a diversionary edges is too large (i.e $\rho \geq 0.6$).

\section{Concluding Remarks}\label{conclusion}
In this paper, we have considered the mathematical modelling, theoretical framework, and computational aspects of recovering a density-based clique relaxation, known as quasi-clique, in networks.
The  planted quasi-clique problem is an extension of the planted clique problem. We have shown that the planted quasi-clique can be solved by relaxing it to a convex programming. This was achieved by borrowing techniques from low rank matrix decomposition and adapting it to the problem under consideration. We have used this formulation to solve the planted maximum quasi-clique problem in the randomized case. In the case where the input graph contains the desired single large  dense subgraph and a moderate number of diversionary vertices and edges, then the relaxation is exact. Therefore, in planted case, these difficult combinatorial optimization problem can be efficiently solved using the tractable relaxation. One of our methodological contributions is the sharp theoretical bounds obtained for the dual matrix. Our proof follows the state of the art in the matrix decomposition literature. However, our innovations lie in the tools we used for our analysis to achieve better results. We improved the results on low rank matrix decomposition by deriving the bound on our dual matrix, using the matrix $l_{\infty, 2}$ norm. This norm has previously been used for matrix completion problem. We established conditions under which recovery is achievable by deriving a dual matrix, certifying the optimality of our solution.  

The results in Section \ref{sec:73} is consistent with our theoretical findings that when $\gamma \geq \frac{c_0 \mu r \log n}{n}$, our convex formulation is guaranteed to recover the planted quasi-clique. However, this is constrained by the level of random noise present. No matter what the value of $\gamma$ is, when $\rho \geq 0.6$ for the cases we considered, recovery is impossible. Further study could consider establishing, theoretically, the upper bound for the amount of random noise that could be added for which recovery is still possible.

\bibliographystyle{plainnat}
\bibliography{main_paper_draft_archived.bib}
\section*{Appendix}
\appendix
\section{Proof of Lemma \ref{lemma2.4candes}}\label{app1}
\begin{proof}
		To prove this lemma, we first show that $(B, C)$ is an optimal pair to the problem (\ref{mdp}) and then show that it is unique. From the optimality condition base on the subgradient condition at $(B, C)$, there must exist a dual matrix, $\mathcal{F}$, which simultaneously belong to the  subdifferential of $B \text{ } (\partial ||B||_*)$ and $C \text{ } (\partial ||C||_1)$. The condition $UU^T + Q = \lambda(Sgn(C) + D)$ of the proposition shows that such matrix exists. Hence $(B , C)$ is an optimal pair. To show uniqueness, we consider a perturbation of $(B, C)$, i.e, $(B + \nu_B , C + \nu _ C)$, which is also a minimizer. Since $B + C = B + \nu_B + C + \nu_C, \nu_B + \nu_c$ must be equal to zero. Now, applying the subgradient condition, we have:
	\begin{align}
	||B + \nu _B||_* + \lambda || C + \nu_C||_1 & \geq ||B||_* + \lambda ||C||_1 + \langle UU^T + Q_0, \nu_B \rangle  \nonumber \\
	& + \lambda \langle Sgn(C) + D_0 , \nu_C \rangle \nonumber \\
	& = ||B||_* + \lambda || C||_1 + \Delta + \delta \label{cross}
	\end{align}
	Recall that $\mathcal{F} \in \partial ||B||_*$ and $\mathcal{F} \in \partial ||C||_1 \implies \mathcal{F} = UU^T + Q = \lambda(Sgn(C) + D)$. Therefore,
	\begin{align*}
	\Delta &= \langle UU^T + Q_0, \nu_B \rangle = \langle UU^T + \mathcal{P}_{T^ \perp}(Q_0), \nu_B \rangle\\
	& = \langle \mathcal{F} - \mathcal{P}_{T^\perp} (Q) + \mathcal{P}_{T^\perp}(Q_0), \nu_B \rangle\\
	& = \langle \mathcal{P}_{T^\perp} (Q_0) - \mathcal{P}_{T^\perp}(Q), \nu_B \rangle +  \langle \mathcal{F} , \nu_B \rangle.
	\end{align*}
	
	Likewise,
	\begin{align*}
	\delta &= \lambda \langle Sgn(C) + D_0, \nu_C \rangle = \langle \lambda Sgn(C) + \lambda \mathcal{P}_{\Gamma^\perp} (D_0), \nu_C \rangle \\
	& = \langle \mathcal{F} - \lambda \mathcal{P}_{\Gamma^\perp}(D) + \lambda \mathcal{P}_{\Gamma ^ C}(D_0), \nu_C \rangle\\
	&= \langle \lambda \mathcal{P}_{\Gamma^\perp} (D_0) - \lambda \mathcal{P}_{\Gamma^\perp}(D), \nu_C \rangle + \langle \mathcal{F}, \nu_C \rangle\\
	&= \lambda \langle \mathcal{P}_{\Gamma^\perp} (D_0) - \mathcal{P}_{\Gamma^\perp} (D), \nu_C \rangle + \langle \mathcal{F}, \nu_C \rangle.
	\end{align*}
	
	Hence,
	\begin{align*}
	\Delta + \delta &= \langle \mathcal{P}_{T^\perp} (Q_0) - \mathcal{P}_{T^\perp}(Q), \nu_B \rangle +  \langle \mathcal{F} , \nu_B \rangle + \lambda \langle \mathcal{P}_{\Gamma^\perp}  (D_0) - \mathcal{P}_{\Gamma^\perp} (D), \nu_C \rangle\\
	& 	+ \langle \mathcal{F}, \nu_C \rangle\\
	&= \langle \mathcal{P}_{T^\perp} (Q_0) - \mathcal{P}_{T^\perp}(Q), \nu_B \rangle + \lambda \langle \mathcal{P}_{\Gamma^\perp}  (D_0) - \mathcal{P}_{\Gamma^\perp} (D), \nu_C \rangle + \langle \mathcal{F}, \nu_B + \nu_C \rangle\\
	&=  \langle \mathcal{P}_{T^\perp} (Q_0) - \mathcal{P}_{T^\perp}(Q), \nu_B \rangle + \lambda \langle \mathcal{P}_{\Gamma^\perp}  (D_0) - \mathcal{P}_{\Gamma^\perp} (D), \nu_C \rangle.
	\end{align*}
	Therefore, equation (\ref{cross}) becomes:
	\begin{align*}
	||B + \nu _B||_* + \lambda || C + \nu_C||_1 & \geq ||B||_* + \lambda ||C||_1 + \langle \mathcal{P}_{T^\perp} (Q_0) - \mathcal{P}_{T^\perp}(Q), \nu_B \rangle + \\
	& \lambda \langle \mathcal{P}_{\Gamma^\perp}  (D_0) - \mathcal{P}_{\Gamma^\perp} (D), \nu_C \rangle  \\
	&= ||B||_* + \lambda ||C||_1 + \left(||\mathcal{P}_{T^\perp}(Q_0)|| - ||\mathcal{P}_{T^\perp}(Q) \right) \left( ||\nu_B ||_* \right) + \\
	& \lambda \left( ||\mathcal{P}_{\Gamma^\perp} (D_0)||_\infty - ||\mathcal{P}_{\Gamma^\perp} (D)||_\infty \right) \left( ||\nu_C||_1 \right)
	\end{align*}
	Since $(Q_0, D_0)$ is any subgradient of $||B||_* + \lambda || C||_1$ at $(B, C)$, we can choose $\mathcal{P}_{T^\perp}(Q_0)$ and $\mathcal{P}_{\Gamma^\perp} (D_0)$ freely inasmuch as they satisfy the following conditions:
	$$||\mathcal{P}_{T^\perp} (Q_0)|| \leq 1, \text{ and } ||\mathcal{P}_{\Gamma^\perp} (D_0)||_\infty \leq 1.$$
	Letting $\mathcal{P}_{\Gamma^\perp} (D_0) = Sgn (\mathcal{P}_{\Gamma^\perp} (\nu_C))$ implies that $||\mathcal{P}_{\Gamma^\perp} (D_0)||_\infty = 1$. Also, since the dual of nuclear norm is the operator norm, there exists a matrix $Q_0$ with $||\mathcal{P}_{T^\perp}(Q_0)|| = 1$ such that $\langle Q_0, \nu_B \rangle = ||\mathcal{P}_{T^\perp} \nu_B||_*$. Therefore, we have
	\begin{align*}
	||B + \nu _B||_* + \lambda || C + \nu_C||_1 & \geq ||B||_* + \lambda ||C||_1 + (1 - ||\mathcal{P}_{T^\perp}(Q)||)(||\mathcal{P}_{T^\perp}(\nu_B)||_*) \\
	&+ \lambda (1 - ||\mathcal{P}_{\Gamma^\perp} (D)||_\infty)(||\mathcal{P}_{\Gamma^\perp} (\nu_C)||_1).
	\end{align*}
	Since $||\mathcal{P}_{T^\perp}(Q)|| < 1$ and $||\mathcal{P}_{\Gamma^\perp} (D)|| < 1$, the last two terms in the equation above are strictly positive except if $||\mathcal{P}_{T^\perp}(\nu_B)||_* = 0$ and $||\mathcal{P}_{\Gamma^\perp} (\nu_C)||_1 = 0$. Hence, $||B + \nu _B||_* + \lambda || C + \nu_C||_1  = ||B||_* + \lambda ||C||_1$ if and only if $\mathcal{P}_{T^\perp}(\nu_B) = 0$ and $\mathcal{P}_{\Gamma^\perp} (\nu_C) = 0$. Note that $\mathcal{P}_{T^\perp}(\nu_B) = \mathcal{P}_{\Gamma^\perp} (\nu_C) = 0$ implies that $\mathcal{P}_T (\nu_B) + \mathcal{P}_\Gamma (\nu _C) = 0$, since $\nu_B + \nu_C = 0$. So that $\mathcal{P}_T(\nu_B) = - \mathcal{P}_\Gamma (\nu_C) = 0$ (since $\Gamma \cap T = \{ 0 \}$). This implies that $\nu_B = \nu_C = 0$. Therefore, $||B + \nu _B||_* + \lambda || C + \nu_C||_1  > ||B||_* + \lambda ||C||_1$ except if $\nu_B = \nu_C = 0$. 
\end{proof}
\end{document}